\documentclass[a4paper,11pt]{amsart}
\usepackage{amsthm,amsmath,amsfonts,amssymb}
\usepackage{color}

\usepackage{cases}
\usepackage[colorlinks=true,urlcolor=teal,
citecolor=magenta,linkcolor=teal,linktocpage,pdfpagelabels,
bookmarksnumbered,bookmarksopen]{hyperref}
\usepackage[capitalize, nameinlink, noabbrev]{cleveref}%DOPO hyperref

\usepackage{graphicx}
\usepackage{subfigure}
\usepackage{sidecap}
\sidecaptionvpos{figure}{t}
\usepackage{overpic}
\usepackage{xfrac}

\usepackage{pgf,tikz}
\usepackage{mathrsfs}
\usetikzlibrary{arrows}

\newtheorem{thm}{Theorem}[section]
\newtheorem{crit}[thm]{Criterion}

\newtheorem{cor}[thm]{Corollary}

\newtheorem{defin}[thm]{Definition}

\newcommand{\R}{{\mathbb{R}}}
\renewcommand{\S}{{\mathbb{S}}}

\newcommand{\Om}{\Omega}

\newcommand{\calH}{\mathcal{H}}

\DeclareMathOperator{\RR}{reach}

\DeclareMathOperator{\dist}{dist}
\DeclareMathOperator{\DIVV}{div}
\newcommand{\de}{\mathrm{d}}
\DeclareMathOperator{\esssup}{ess\, sup}

\begin{document}

\definecolor{ffffff}{rgb}{1.,1.,1.}
\definecolor{cqcqcq}{rgb}{0.75,0.75,0.75}
\definecolor{uuuuuu}{rgb}{0.25,0.25,0.25}

\title[Geometric criteria $\&$ existence of capillary surfaces]{Geometric criteria for the existence of capillary surfaces in tubes}
\author[G.~Saracco]{Giorgio Saracco}
\address[Giorgio Saracco]{Dipartimento di Matematica e Informatica ``Ulisse Dini'' (DIMAI), Universit\`a di Firenze, viale Morgagni, 67/A, 50134 Firenze (FI), Italy}
\email{giorgio.saracco@unifi.it}

\dedicatory{In memoriam of Robert Finn}

 %53A10, 76D45, 58E12


\begin{abstract}
We review some geometric criteria and prove a refined version, that yield existence of capillary surfaces in tubes $\Omega\times \R$ in a gravity free environment, in the case of physical interest, that is, for bounded, open, and simply connected $\Omega \subset \R^2$. These criteria rely on suitable weak one-sided bounds on the curvature of the boundary of the cross-section $\Omega$.
\end{abstract}

 \hspace{-2cm}
 {
 \begin{minipage}[t]{0.6\linewidth}
 \begin{scriptsize}
 \vspace{-3cm}
 This is a pre-print of an article published in \emph{Expo.\ Math.}. The final authenticated version is available online at: \href{https://doi.org/10.1016/j.exmath.2024.125547}{https://doi.org/10.1016/j.exmath.2024.125547}
 \end{scriptsize}
\end{minipage} 
}

\maketitle

\section{Introduction}

\begin{SCfigure}[60]
\begin{overpic}[height=8cm]{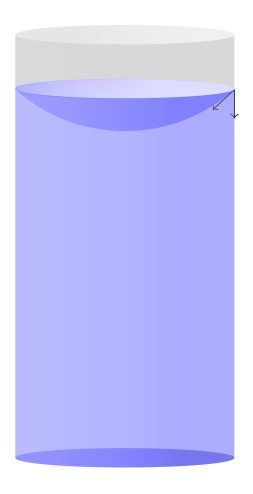}
\put(43,75){$\gamma$}
\put(13,79.5){$\Gamma$}
\end{overpic}
\caption{A fluid poured in a thin tube gives rise to a capillary surface $\Gamma$, that is, a (graphical) surface with constant mean curvature adhering to the walls of the tube with a constant contact angle $\gamma$ determined by the properties of the fluid, the air, and the material of the tube.}
\label{fig:capillary_tube}
\end{SCfigure}

In this brief note, we are interested in some purely geometrical criteria that allow to determine whether capillary surfaces exist in a cylindrical tube of cross-section a bounded, open, and simply connected subset $\Omega$ of $\R^2$. Let us suppose that in the cylinder $\Omega \times \R$ (closing one of the ends) there are two immiscible and incompressible phases in equilibrium (e.g., air and water) separated by an interface $\Gamma$, and let us assume that this one can be represented by the graph of a function $u$. Then, the energy of this physical system, whenever gravity is absent or can be neglected (for instance, on Earth's surface whenever the diameter of the cross-section is sufficiently small and the mass of the fluid is small), consists of the sum of three different terms: the free surface energy that represents the work necessary to build the separating interface; the wetting energy that quantifies the work of the adhesion forces between the phases and the rigid vertical walls of the cylinder; a volume constraint standing for the finiteness of the mass of the fluid we are considering. Hence, up to the multiplicative factor of the surface tension, the energy of the system is
\begin{equation}\label{eq:energy}
\int_\Omega \sqrt{1+|\nabla u|^2}\,\de x - \cos(\gamma) \int_{\partial \Omega} u\, \de \calH^{1}(x) +\int_\Omega \lambda u\,\de x \,,
\end{equation}
where the first term is the surface energy, the second one the adhesion energy and where $\gamma$ is the contact angle measured inside the lower fluid between the phases and the cylinder, and the last one represents the volume constraint, being $\lambda$ a Lagrange multiplier. The energy expression~\eqref{eq:energy} had already been derived by Gauss, unifying previous theories by Young and Laplace. A nice account and modern derivation is available in~\cite{Fin01}, and in~\cref{fig:capillary_tube} there is a sketch of the physical situation we are interested in. 
Writing down the Euler--Lagrange equation of the functional~\eqref{eq:energy}, one finds that smooth critical points need to satisfy
\begin{alignat}{4}
%\begin{cases}
\DIVV (Tu) &= \lambda \,, \qquad \qquad &&\text{in $\Om$}, \label{eq:cap} \\
Tu \cdot \nu_\Om  &= \cos(\gamma)\,, &&\text{on $\partial\Om$}\,, \label{eq:bc}
%\end{cases}
\end{alignat}
where $Tu$ is the vector field
\[
Tu = \frac{\nabla u}{\sqrt{1+|\nabla u|^2}}\,,
\]
and $\nu_\Omega$ the outward normal to $\Omega$. A caveat is that when stating the above mathematical problem, one considers a cylinder of infinite length, that is an unrealistic physical situation. Nevertheless, if the mathematical formalism leads to (and it does lead to) a solution that is the graph of a function $u$ bounded from below, then a physically meaningful solution can be derived by adding a suitable constant, that is, by closing one end of the cylinder and covering it with fluid.

The lone PDE~\eqref{eq:cap} without the boundary condition~\eqref{eq:bc} is generally referred to as the \emph{prescribed mean curvature equation} because the term in the LHS, $\DIVV(Tu)$, represents the pointwise mean curvature of the graph of $u$. Hence, capillary surfaces have constant mean curvature given by $\lambda$. For the sake of explanation, let us momentarily disregard the boundary condition~\eqref{eq:bc} and let us consider the $1$-dimensional case, that is, $\Omega$ is, without loss of generality, the interval $(-a,a)$. When searching for a $\mathrm{C}^2$ solution, we are looking for a continuously twice-differentiable function $u$ such that the curve $(x,u(x))$ has constant curvature given by $\lambda$. Since in dimension $1$ the only curves with (positive) constant curvature $\lambda$ are arcs of circles of radius $\lambda^{-1}$, our solution has to be one of such arcs. This imposes some geometric restraint on the initial choice of $\lambda$. Indeed, in order to be able to bridge the gap spanned by the interval $(-a,a)$ with an arc of a circle of radius $\lambda^{-1}$, we necessarily need $\lambda^{-1}$ to be \emph{at least} $a$, refer also to \cref{fig:bridging}. 
\begin{figure}[t]%
\centering
	\subfigure[Bridging a gap wide $2a$ via an arc of circle requires a radius of at least $a$. \label{fig:bridging}]
	{
	\begin{tikzpicture}
	\draw[black] (-2,0) -- (2,0);
	\draw[black, dashed] (-2,0) -- (-2,4);
	\draw[black, dashed] (2,0) -- (2,4);
	\draw (-2.2,.3) node {$(-a,0)$};
	\draw (1.95,.3) node {$(a,0)$};
	\clip (-2,0) rectangle (2,4);
	\draw[black] (.7,4) circle [radius=3cm];
	\end{tikzpicture}
	}
	\hspace{.5cm}
	\subfigure[Prescribing the boundary condition~\eqref{eq:bc} forces the symmetry of the solution. \label{fig:bridging_sym}]
	{
	\begin{tikzpicture}
	\draw[black] (-2,0) -- (2,0);
	\draw[black, dashed] (-2,0) -- (-2,4);
	\draw[black, dashed] (2,0) -- (2,4);
	\draw (-2.2,.3) node {$(-a,0)$};
	\draw (1.95,.3) node {$(a,0)$};
	\clip (-2,0) rectangle (2,4);
	\draw[black] (0,4) circle [radius=3cm];
	\end{tikzpicture}
	}
\caption{Two solutions of the prescribed curvature equation in $1$d. On the left one without enforcing any boundary condition, while on the right with a Neumann-like condition.}
\label{fig:bridges}
\end{figure}
Hence, a necessary condition for existence is that
\begin{equation}\label{eq:1-d}
\lambda \le \frac 1a = \frac{\mathcal{H}^0((-a,a))}{\mathcal{H}^1((-a,a))}.
\end{equation}
This is not a condition peculiar to dimension $1$. In the physical situation\footnote{The following holds in general dimension $N$, using the relevant Hausdorff measures.} of a bounded, open, and simply connected $\Omega\subset \mathbb{R}^2$, integrating the PDE~\eqref{eq:cap} on $\Omega$ and using the Gauss--Green Theorem (assuming $\Omega$ Lipschitz), one has
\[
-\int_{\partial\Omega} \frac{\nabla u \cdot \nu_\Omega}{\sqrt{1+|\nabla u|^2}}\, \de \mathcal{H}^{1}(x) = \int_\Omega \DIVV(Tu)\, \de x = \lambda\int_\Omega 1\,\de x = \lambda|\Omega|.
\]
Taking now the absolute value on the LHS, moving it under the integral sign and using that $|Tu|\le 1$, one finds the necessary condition to existence 
\begin{equation}\label{eq:necessary_no_bc}
\lambda \le \frac{P(\Omega)}{|\Omega|}\,,
\end{equation}
where $P(\,\cdot\,)$ denotes the distributional perimeter, which for Lipschitz sets $E$ coincides with $\mathcal{H}^1(\partial E)$, and we refer to~\cite{AFP00book, Mag12book} for the basic facts of sets of finite perimeter.
Further, carrying out the same steps on any Lipschitz subset $E$ compactly contained in $\Omega$ one also gets as necessary condition that
\begin{equation}\label{eq:necessary_compact_subsets}
\lambda <  \frac{P(E)}{|E|}\,,
\end{equation}
where the strict inequality comes from the fact that the vector field $Tu$ is such that $|Tu|<1$ on any subset compactly contained in $\Omega$.
 
Let us now consider the PDE~\eqref{eq:cap} coupled with the boundary condition~\eqref{eq:bc}. In this case we do not get an upper bound on $\lambda$ like in~\eqref{eq:necessary_no_bc}--\eqref{eq:necessary_compact_subsets}, rather it is implicitly determined by the geometrical nature of the problem. For the sake of simplicity, we start again considering the $1$-dimensional case. Requiring~\eqref{eq:bc} on the boundary forces the symmetry of the arc with respect to the axis of the segment $(-a,a)\times\{0\}$, see also \cref{fig:bridging_sym}. Thus, such arc belongs to a circle centered on the $y$-axis at $(0,y_C)$ and, by symmetry, it intersects the walls of the cylinder at $(-a,k)$ and $(a,k)$, for some $k$. Without loss of generality, one can assume the height $k$ to be fixed. It is easy to see that the function that to each $y$ associates the angle created by the circle centered at $(0,y)$ with the line $\{a\}\times \mathbb{R}$ at $(a,k)$ is strictly monotonic, therefore the angle $\gamma$ prescribed by~\eqref{eq:bc}, coupled with any choice of the height $k$ completely determines the height $y_C=y_C(k)$ of the center. Changing the initial height $k$ modifies $y_C$ but not their relative distance. Hence, the radius of the circle ends up being determined by the geometry of the problem, and in turns so it is the prescribed curvature $\lambda$. The only degree of freedom that one has is the choice of the height $k$, which mathematically corresponds to the uniqueness of the solution up to a vertical translation.

The same happens in the higher dimensional case. Integrating the PDE~\eqref{eq:cap} on $\Omega$, owing to Gauss--Green Theorem and taking into account the boundary condition, in place of the large inequality~\eqref{eq:necessary_no_bc}, one gets that $\lambda$ needs to satisfy
\begin{equation}\label{eq:critical}
\lambda = \cos(\gamma) \frac{P(\Omega)}{|\Omega|}\,.
\end{equation} 
Performing the same reasoning on any proper Lipschitz subset $E$ of $\Omega$, taking into account the boundary condition~\eqref{eq:bc} on $\partial E \cap \partial \Omega$ and that $|Tu|<1$ on $\partial E\cap \Omega$, one finds the following necessary condition to existence 
\begin{equation}\label{eq:necessary}
\lambda <  \frac{P(E; \Omega) +\cos(\gamma) P(E; \partial \Omega)}{|E|}\,,
\end{equation}
which is the analog of~\eqref{eq:necessary_compact_subsets}, and where $P(E; A)$ denotes perimeter of $E$ relative to the set $A$. 
Pairing~\eqref{eq:critical} with~\eqref{eq:necessary} one finds a necessary condition that is purely geometrical, that is, the strict inequality
\begin{equation}\label{eq:necessary_general_gamma}
\cos(\gamma) \frac{P(\Omega)}{|\Omega|} < \frac{P(E; \Omega) +\cos(\gamma) P(E; \partial \Omega)}{|E|}
\end{equation}
must hold for all Lipschitz proper subsets $E$ of $\Omega$. Thus, existence of capillary surfaces is intimately tied to a purely geometrical problem set in one lower dimension, which has been studied by Concus and Finn in~\cite{CF84, Fin84}, refer also to~\cite[Chap.~19]{Mag12book}. It is easy to see that if condition~\eqref{eq:necessary_general_gamma} is met for a vertical contact angle, that is $\gamma=0$, then it also is for any $\gamma \in [0, \sfrac \pi2]$. If $\gamma=\sfrac \pi 2$ it is trivial; if otherwise $\gamma\in [0,\sfrac \pi 2)$ we can divide both terms in~\eqref{eq:necessary_general_gamma} by $\cos(\gamma)$ and it suffices to notice that the RHS is strictly increasing as a function of $\gamma$.

For this reason, we focus on the choice $\gamma=0$, in which case the necessary condition becomes
\begin{equation}\label{eq:geo_necessary_1}
\frac{P(\Omega)}{|\Omega|} < \frac{P(E)}{|E|}\,,
\end{equation}
for all Lipschitz proper subsets $E$ of $\Omega$. Such a condition was first shown to be necessary by Concus and Finn~\cite[Lem.~4]{CF74}, under the Lipschitz regularity of $\Omega$, required to employ the Gauss--Green Theorem. Via standard arguments, this can be weakened to a piecewise Lipschitz request on $\Omega$. Under this regularity assumption, inequality~\eqref{eq:geo_necessary_1} was later shown in Giusti's seminal paper~\cite{Giu78} to be necessary for all choices of proper subsets $E\subsetneq \Omega$ of locally finite perimeter, considerably weakening the Lipschitz request on the subsets $E$ via an approximation argument. In the same paper Giusti proves the sufficiency of such a condition, see also the comprehensive treatise~\cite[Chap.~6]{Fin86book}. The intrinsic geometric nature of this condition led the author to ask himself whether some geometric criterion could be proved to ensure the validity of~\eqref{eq:geo_necessary_1}, in~\cite[Thm.~A.1 and Cor.~A.1]{Giu78} he proved one, under the assumptions that $\Omega$ is bounded, open, and convex, and of class $\mathrm{C}^1$, and this will be touched upon in \cref{sec:Giusti}. We remark that this first criterion provides an \emph{``if and only if''} statement. Interestingly, one cannot just drop the convexity hypothesis or replace it with star-shapedness and hope to retain the same criterion, as shown by means of counterexamples by Finn and Giusti~\cite{FG79b}. Few years later Chen~\cite[Thm.~4.1]{Che80} was able to extend the criterion to bounded, open, simply connected, and piecewise Lipschitz sets that enjoy a ``strict'' interior ball condition, where the meaning of ``strict'' will be made clear later on in \cref{sec:Chen}, where we discuss this extension. Notably, dropping convexity and replacing it with this weaker assumption, produces a sufficient but not anymore necessary criterion, unless an extra assumption is requested. In the last decade new Gauss--Green formulas for much less regular sets $\Omega$ have been proved, we here specifically refer to~\cite{LS20a}, but the topic has been very active and a far from complete list is~\cite{BCM21, CCT19, CTZ09, CdC18, CdC19}. With these new formulas at one's disposal, it has been possible to prove existence of solutions of~\eqref{eq:cap}--\eqref{eq:bc} under much weaker regularity conditions on $\Omega$, see~\cite{LS18a}, covering also wildly irregular sets, for instance balls with infinitely many holes accumulating toward their boundaries~\cite[Sect.~3]{LS18b}. Namely, researchers went from asking the piecewise Lipschitz regularity of $\Omega$ to the request
\begin{equation}\label{eq:P=H}
P(\Omega) = \calH^1(\partial \Omega)\,,
\end{equation}
and the existence of a positive constant $k$ depending only on $\Omega$ such that for all subsets $E\subset \Omega$ one has
\begin{equation}\label{eq:Poincare}
\min\{\,P(E; \partial \Omega); P(\Omega \setminus E; \partial \Omega) \,\} \le k P(E; \Omega)\,,
\end{equation}
that is, a Poincar\'e-type inequality. Whenever $\Omega$ is such that it satisfies~\eqref{eq:P=H} and~\eqref{eq:Poincare}, inequality~\eqref{eq:geo_necessary_1} again provides a necessary and sufficient condition to existence, see~\cite[Thms.~4.3 and~4.7]{LS18a}. It is in this new framework that we shall see that, combining results from~\cite{Sar18, Sar21}, a sufficient criterion for existence of solutions similar to that of Chen~\cite{Che80} can be proved in a wider generality, dropping the piecewise Lipschitz requests on the boundary of $\Omega$, and only assuming a suitable weak condition on the curvature (in terms of the \emph{reach}), and this will be the topic of \cref{sec:Saracco}.

The plan of the paper is the following. In \cref{sec:Cheeger} we introduce a related, useful problem through which we can state the problem in a ``different language''. In \cref{sec:Giusti} and \cref{sec:Chen} we review the first two geometric criteria yielding existence of capillary surfaces. In \cref{sec:Saracco} we exploit few recent results to prove a refined version of these criteria that has been hinted to in~\cite{Sar21} but neither formally stated nor proved.

\section{A related problem}\label{sec:Cheeger}

Let us set aside for a moment the capillarity problem~\eqref{eq:cap}--\eqref{eq:bc}, with all the physical implications it carries along, and let us consider a related problem: the lone prescribed mean curvature differential equation without any boundary datum
\begin{equation}\label{eq:pmc}
\DIVV (Tu) = H \,, \qquad \text{in $\Om$}\,,
\end{equation}
where $H$ is a fixed positive constant. Reasoning as we did in the introduction, by an application of the Gauss--Green Theorem, and using the approximation argument of Giusti~\cite[Sect.~1]{Giu78}, we have as necessary condition to existence that
\begin{equation}\label{eq:nec_pmc}
H < \frac{P(E)}{|E|}\,,
\end{equation}
for all proper subsets $E\subset \Omega$ of locally finite perimeter. Then, one has two cases: either the strict inequality in~\eqref{eq:nec_pmc} holds as well for $\Omega$ in place of $E$ (non-critical case); or the curvature $H$ equals the ratio perimeter over volume of $\Omega$ (critical case). In both cases, condition~\eqref{eq:nec_pmc} can again be proved to be sufficient for existence of solutions when $\Omega$ is bounded, open, simply connected, and either it is piecewise Lipschitz or it satisfies~\eqref{eq:P=H}--\eqref{eq:Poincare}.

The interesting bit is that in the critical case any solution of~\eqref{eq:pmc} will be a capillary surface, since it automatically satisfies~\eqref{eq:cap}--\eqref{eq:bc} with a vertical contact angle, that is, for $\gamma=0$ (the boundary datum being assumed almost everywhere). Even more, such a solution is unique up to vertical translations. For the sake of convenience, we sum up this result in the next statement (refer to~\cite[Thm.~2.1]{Giu78} for the piecewise Lipschitz case and to~\cite[Thm.~5.1]{LS18a} for the general case combined with~\cite[Sect.~5.1]{LS18b}).

\begin{thm}\label{thm:pmc-iff}
Let $\Omega$ be a bounded, open, and simply connected subset of $\mathbb{R}^2$ satisfying~\eqref{eq:P=H}--\eqref{eq:Poincare}, and let $H\in \mathbb{R}$ fixed. The following are equivalent:
\begin{itemize}
\item[(a)] \eqref{eq:nec_pmc} holds and $H= P(\Omega)|\Omega|^{-1}$;
\item[(b)] there exists a unique (up to translations) solution of~\eqref{eq:pmc};
\item[(c)] there exists a solution $u$ such that $Tu\cdot \nu_\Omega =1$ a.e.\ on $\partial \Omega$, i.e., that solves~\eqref{eq:cap}--\eqref{eq:bc} with $\gamma=0$.
\end{itemize}
\end{thm}

Since~\eqref{eq:nec_pmc} needs to hold for all subsets of $\Omega$ of locally finite perimeter, if we take the infimum in~\eqref{eq:nec_pmc} among all Borel subsets of $\Omega$, we find
\[
H \le \inf \left\{\,\frac{P(E)}{|E|}\,:\, |E|>0 \,\right\} = h(\Omega)\,,
\]
where the constant on the RHS is known as the \emph{Cheeger constant of $\Omega$}, and it was first introduced in~\cite{Che70}. The problem of computing the constant $h(\Omega)$ and of determining the subsets of $\Omega$ attaining it (called \emph{Cheeger sets of $\Omega$}) has gained a lot of attention in the last decades. Whenever $\Omega$ is the unique Cheeger set in itself, we shall call it a \emph{minimal Cheeger set}. We refer to the two surveys~\cite{Leo15, Par11}, and to~\cite{KL06, LNS17, LP16, LS20b, Sar21} for results in the 2-dimensional case that we are interested in, and the references therein.

Introducing the constant and being aware of the results on it is very useful, since we can easily restate the non-critical and the critical case in terms of the Cheeger constant and of Cheeger sets. Indeed, whenever $H< h(\Omega)$ we are in the non-critical case, while when $H=h(\Omega)$ \emph{and} $h(\Omega)$ is uniquely attained by $\Omega$, that is, when $\Omega$ is a minimal Cheeger set, we are in the critical case.

In view of this parallel, and the equivalence of the critical case for~\eqref{eq:pmc} (\cref{thm:pmc-iff}~(a)) with the capillarity problem~\eqref{eq:cap}--\eqref{eq:bc} with $\gamma=0$ (\cref{thm:pmc-iff}~(c)), what one needs to find are criteria on a bounded, open, simply connected, and planar set $\Omega$ (either piecewise Lipschitz, or satisfying~\eqref{eq:P=H}--\eqref{eq:Poincare}) that ensure that $\Omega$ is a minimal Cheeger set. The proofs we will provide of the criteria contained in the next sections exploit this deep link between the two problems, and we shall see how criteria for $\Omega$ being a minimal Cheeger set almost immediately port to criteria for existence of solutions of~\eqref{eq:cap}--\eqref{eq:bc} with $\gamma=0$.

We also mention that such self-minimality criteria for the Cheeger problem appear to be useful in other contexts such as image reconstruction problems as the ROF model~\cite{ROF92}, failure of planar plates subject to a vertical load~\cite{Kel80}, and viscoplastic fluids~\cite{HG23, HG23s}.

For the sake of completeness, we remark that minimal Cheeger sets are in some contexts called \emph{calibrable sets}, see~\cite[Def.~3 and Rem.~6]{BCN02}. We also mention that the task of verifying the necessary condition~\eqref{eq:geo_necessary_1}, or, as discussed, determining whether $\Omega$ is a minimal Cheeger set, is equivalent to finding a vector field on $\Omega$ with some special properties, and we refer the interested reader to~\cite[Thms.~1 and~2]{Fin79}, and also to~\cite[Thm.~3]{Gri06} and~\cite{Str08, Str10}.

\section{Convex sets}\label{sec:Giusti}

The first criterion for existence that appeared dates back to 1978 and it is due to Giusti~\cite[Thm.~A.1 and Cor.~A.1]{Giu78}. Several years later two criteria in different settings were independently proved, and these paired with the observations of \cref{sec:Cheeger} immediately allow to recover Giusti's one. In particular, we refer to Bellettini--Caselles--Chambolle~\cite[Rem.~6 and Thm.~4]{BCN02} in the framework of $\mathrm{C}^{1,1}$, calibrable, and convex sets, and to Kawohl--Lachand-Robert~\cite[Thm.~2]{KL06} in the language of Cheeger sets. This latter, paired with the observations of \cref{sec:Cheeger}, provides a very quick proof of the criterion.

\begin{crit}\label{crit:Giusti}
Let $\Omega \subset \R^2$ be a bounded, open, and convex set, and let $\bar \kappa$ be defined as
\[
\bar \kappa = \esssup \kappa_{\partial \Omega}\,,
\]
where $\kappa_{\partial \Omega}$ represents the curvature of $\partial \Omega$. Then the PDE~\eqref{eq:cap} with boundary condition~\eqref{eq:bc} has a solution for $\gamma =0$ \emph{if and only if} 
\begin{equation}\label{eq:curvature_condition}
\bar \kappa \le \frac{P(\Omega)}{|\Omega|}\,.
\end{equation}
\end{crit}

\begin{proof}
By~\cite[Thm.~2]{KL06} a bounded, open, and convex planar set $\Omega$ is self-Cheeger if and only if $\bar \kappa \le P(\Omega)/|\Omega|$. Moreover, convex sets have a unique Cheeger set~\cite{AC09}, hence it also tells us that $h(\Omega)$ is \emph{uniquely} attained by $\Omega$. 

In turns, this says we are in the critical case for solving~\eqref{eq:pmc}, that is, statement~(a) in \cref{thm:pmc-iff} holds. Since convex implies Lipschitz, the regularity assumptions on $\Omega$ requested in \cref{thm:pmc-iff} are satisfied. Therefore, by the equivalency stated by the theorem, there is a solution of~\eqref{eq:cap}--\eqref{eq:bc}, with $\gamma =0$ (\cref{thm:pmc-iff}~(c)).
\end{proof}

It is here useful to define what we mean by curvature of a convex set, which a priori is only Lipschitz, and thus the curvature may not be defined in the classical sense, that requires a $\mathrm{C}^2$ regularity. 

Given a convex set $\Omega$, its \emph{support function} is $p : \S^1 \to \R$ defined as
\[
p(\theta) = \sup_{(x,y)\in \Omega} \{\, x\cos(\theta) + y\sin(\theta)\,\},
\]
which is Lipschitz continuous and allows to identify the hyperplane orthogonal to $(\cos(\theta), \sin(\theta))$ supporting the convex set $\Omega$. The boundary of $\Omega$ can be then described almost everywhere as the simple and closed curve $(x(\theta), y(\theta))$, with $\theta\in [0, 2\pi]$, given by\begin{align*}
x(\theta) &= p(\theta)\cos(\theta) -p'(\theta)\sin(\theta),\\
y(\theta) &= p(\theta)\sin(\theta) + p'(\theta)\cos(\theta).
\end{align*}
If a convex set is of class $\mathrm{C}^2$, its support function $p$ is twice differentiable and the curvature radius $\rho$ is such that $\rho = p+p''>0$, and the curvature is its reciprocal. Conversely, given a Radon measure $p : \S^1 \to \R$, satisfying $p+ p'' \ge 0$ in a distributional sense, one can find a convex set $\Omega$, whose support function is given by $p$. One can prove that there is a bijective correspondence between convex sets and Radon measures $p$ on $\S^1$ such that $p+p''\ge 0$, and the curvature $\kappa$ can be defined in general as the ratio $1/\rho$, being $\rho = p+p''$. The supremum $\bar \kappa$ is then defined as $+\infty$ if $\kappa$ is not bounded, otherwise as the supremum of the Lebesgue precise representative of $\kappa$.

What is important to note is that if $\bar \kappa$ is finite, then the classical fact for $\mathrm{C}^2$ convex sets $\Omega$ that at any point $x\in \partial \Omega$ there exists a ball through $x$ of radius $1/\bar \kappa$ and interior to $\Omega$ remains true. This in particular allows to restate \cref{crit:Giusti} in terms of an interior ball condition of radius $P(\Omega)/|\Omega|$ at any point of the boundary $x\in \partial \Omega$. The criteria we shall discuss in the next sections are indeed essentially stated through this property.

For the sake of completeness, we note that such a criterion has been proved also in different frameworks, when the underlying metric is non Euclidean, but rather anisotropic. In particular, we refer the interested reader to~\cite[Cor.~5.3]{KN08}, see also~\cite[Thm.~8.1]{BNP01}. These criteria have been later extended for general $N$-dimensional convex sets, see~\cite{ACC05} for the Euclidean case and~\cite{CCMN08} for the anisotropic case.

\section{Non convex piecewise Lipschitz sets}\label{sec:Chen}

On the one hand, it was almost immediately clear that the lone curvature condition~\eqref{eq:curvature_condition} is not enough to guarantee the existence of solutions, and that in the previous statement convexity played a major role. Indeed, in~\cite[Sec.~1]{FG79b}, the authors show that existence may fail when negative curvatures occur, even if they are small compared to the ratio $P(\Omega)/|\Omega|$, and even if one requires starshapedness. Specifically, Finn and Giusti consider two balls of radii $R$ and $r$, with $R>r$, with non empty intersection, as shown in \cref{fig:counterexample}.
\begin{figure}
\centering
\begin{tikzpicture}
\begin{scope}[scale=2]
\node at (-1,.75) {$\Omega$};
\draw[thick](0,0) circle (1cm);
\draw[thick](1.5,0) circle (.75cm);
\fill[white] (0,0) circle (1cm);
\fill[white] (1.5,0) circle (.75cm);
\draw[dash pattern=on 1pt off 1pt, thin] (0,0) -- (0.895, 0.444);
\node at (0.35,0.35) {$R$};
\draw[dash pattern=on 1pt off 1pt, thin] (1.5,0) -- (0.895, 0.444);
\node at (1.2,0.35) {$r$};
\draw[dash pattern=on 1pt off 1pt, thin] (0,0) -- (1.5, 0);
\draw[dash pattern=on 1pt off 1pt] (0.2, 0) arc (0:27:0.2);
\node at (0.45,0.1) {$\theta_R$};
\draw[dash pattern=on 1pt off 1pt] (1.3, 0) arc (180:143:0.2);
\node at (1.2,0.1) {$\theta_r$};
\end{scope}
\end{tikzpicture}
\caption{A set $\Omega$ that, when smoothed out, satisfies the curvature condition~\eqref{eq:curvature_condition}, yet existence of solutions of the PDE~\eqref{eq:cap} with boundary condition~\eqref{eq:bc} for $\gamma=0$ fails.}\label{fig:counterexample}
\end{figure}
One can smooth out the intersection points in such a way that the curvature $\kappa$ is strictly less than $\sfrac 1r$ at all points of the boundary. Calling $\theta_R$ and $\theta_r$ the angles drawn by the segment through the two centers and by those through the centers and the intersection point (refer to \cref{fig:counterexample}), and disregarding higher order terms in $\theta_R$ and $\theta_r$, one has
\begin{equation}\label{eq:rapRr}
\frac{P(\Omega)}{|\Omega|} \approx \frac{2\pi(R+r)-2(\theta_R R + \theta_r r)}{\pi(R^2 + r^2)} \to 2\frac{R+r}{R^2+r^2}\,,
\end{equation}
for sufficiently small angles (that is, sufficiently spaced far apart centers). Thus, the curvature condition~\eqref{eq:curvature_condition} is satisfied, provided that
\[
\frac 1r < 2 \frac{R+r}{R^2+r^2}\,,
\]
that is, if $r>R(\sqrt{2}-1)$. Nevertheless, for such values the necessary condition~\eqref{eq:geo_necessary_1} for existence fails to hold, since the ball $B_R\subset \Omega$ has a strictly better perimeter to volume ratio than $\Omega$ itself. Indeed, in view of~\eqref{eq:rapRr},
\[
\frac{P(B_R)}{|B_R|} = \frac 2R < 2\frac{R+r}{R^2+r^2} \approx \frac{P(\Omega)}{|\Omega|}\,,
\]
provided that the centers of the two balls are sufficiently far apart.

On the other hand, also the curvature plays a special role, as implied by~\cite[Thm.~3]{FK02} that states that, for $\Omega$ sufficiently smooth, if there exists even a single point on the boundary where the curvature is strictly greater than $P(\Omega)/|\Omega|$, then the necessary condition~\eqref{eq:geo_necessary_1} fails.

What Chen realized is that the key point was neither the bound on the curvature nor the convexity, rather what the two condition paired together implied: the existence of an interior ball of radius $|\Omega|/P(\Omega)$ through any point of the boundary. As he allowed also less regular sets, that is, piecewise Lipschitz, he stated the condition as an interior ``rolling'' ball condition to take into account that at non regular points $x$ of $\partial \Omega$ one has a cone of inward normals, in place of a uniquely defined one. In particular, since the boundary is assumed to be piecewise Lipschitz, one can easily define a leftmost inward normal $\nu^-$ (resp., a rightmost one $\nu^+$) as the limit of the normals approaching from the left (resp., the right); the resulting cone shall consist of all directions inbetween $\nu^-$ and $\nu^+$ locally pointing toward the inside of the set $\Omega$.

We can now recollect this loose idea into a definition, based on the original one~\cite[Def.~4.1]{Che80}.
\begin{defin}
Let $\Omega\subset \R^2$ be a bounded, open, simply connected, and piecewise Lipschitz set. We say that it enjoys the interior rolling ball condition of radius $r$, if for any $x\in \partial \Omega$, any $\nu \in \S^{1}$ in the cone of interior normals to $\Omega$ at $x$, one has $B_r(x+r \nu) \subset \Omega$. We say that it enjoys the strict interior rolling ball condition if additionally no pair of antipodal points in $\partial B_r(x+r\nu)$ belongs to $\partial \Omega$.
\end{defin}
Loosely speaking, the definition above means that one can ``roll''---hence the adjective ``rolling''---along the boundary of $\Omega$ and internally to $\Omega$ a ball of the given radius $r$, and this can be thought of as a one-sided bound on the curvature of $\partial \Omega$. For the sake of completeness, we remark that Chen only gave the ``strict'' definition without naming it so; the distinction between ``strict'' and ``non strict'' came much later in~\cite{Sar21}. Chen proved that the strict condition, for $r=|\Omega|/P(\Omega)$, is sufficient for existence of solutions, albeit not necessary. In particular, the union of two overlapping balls with the same radius, that is, a situation like that of \cref{fig:counterexample} with $r=R$, with the centers suitably spaced far apart, provides a set for which this condition is not met, while the necessary condition~\eqref{eq:geo_necessary_1} holds, and thus existence is nevertheless ensured, see~\cite[Ex.~5.3]{Che80}.

Let us notice how the interior rolling ball condition of radius $r$ implies that the inner parallel set at distance $r$, that is,
\[
\Omega^r = \{\,x\in \Omega\,:\, \dist(x; \partial \Omega)\ge r \,\}\,,
\]
is simply connected (because $\Omega$ is) and path connected. This path connectedness property turns out to be the key to make the criterion necessary. Let us lay down some more jargon.
\begin{defin}
Let $\Omega\subset \R^2$ be a bounded, piecewise Lipschitz, and simply connected set. We say that it has no necks of radius $r$ if $\Omega^r$ is path connected.
\end{defin}
Chen's definition, refer to~\cite[Def.~5.1]{Che80}, is different from the one given here, but completely equivalent. The one we provided has its roots in~\cite[Def.~1.2 and Rem.~1.3]{LNS17}, where we drew inspiration from the original one. We are now ready to state Chen's criterion~\cite[Thms.~4.1 and~5.2]{Che80}, whose proof we omit, as it is implied by the criterion we state and prove in the next section. 

\begin{crit}\label{crit:Chen}
Let $\Omega\subset \R^2$ be a bounded, piecewise Lipschitz, and simply connected set. 

\begin{itemize}
\item[(i)] If it enjoys the strict interior rolling ball condition for $r=|\Omega|/P(\Omega)$, then the PDE~\eqref{eq:cap} with boundary condition~\eqref{eq:bc} has a solution for $\gamma =0$. 

\item[(ii)] If $\Omega$ has no necks of radius $r=|\Omega|/P(\Omega)$, then the PDE~\eqref{eq:cap} with boundary condition~\eqref{eq:bc} has a solution for $\gamma =0$ \emph{if and only if}  $\Omega$ enjoys the strict interior rolling ball condition for $r=|\Omega|/P(\Omega)$.
\end{itemize}
\end{crit}

We remark that one cannot replace the strict condition with the non strict one, as otherwise the statement would not hold. This is highlighted by the ``Pinocchio'' example shown in \cref{fig:pinocchio}, also called ``keyhole'' or ``proboscis'' by Finn in~\cite{FL94, FM96}, where he studied similar configurations for $\gamma\neq 0$. Such a set is given by the union of the ball $B_1$ centered at the origin and a ball of radius $r=\sin\theta$ centered at $(\cos\theta, 0)$. There exists a choice of $\theta\in (0, \sfrac \pi2)$, that we denote by $\theta_0$ and the corresponding radius by $r_0$ (roughly, $\theta_0\approx 0.531$), such that this set is a minimal Cheeger set. Then for any $T>0$, it is easy to see that the Pinocchio set
\[
\mathcal{P}_T = B_1 \cup \bigcup_{\tau \in [0, T]} B_{r_0}(\cos \theta_0 + \tau, 0)
\]
satisfies the (non strict!) interior rolling ball condition for $r=|\mathcal{P}_T|/P(\mathcal{P}_T)$, it is a Cheeger set in itself, but not a minimal one, since for any $t\in [0, T)$ the proper subset $\mathcal{P}_t$ is such that
\[
\frac{P(\mathcal{P}_t)}{|\mathcal{P}_t|} = \frac{P(\mathcal{P}_T)}{|\mathcal{P}_T|}\,.
\] 
The full computations are available in~\cite[Ex.~4.6]{LP16}, but such an example appeared several times, see~\cite[Sects.~6.13,~6.14, and~6.17]{Fin86book}, and~\cite[Sects.~4.2,~4.3, and~4.4]{BP96}. In particular, one can also smooth out the set controlling the curvature, and build in this way a counterexample to \cref{crit:Giusti} when convexity is removed, in the same spirit of the one described at the beginning of this section, see~\cite[Exs.~6.5 and~6.6]{ABT15}.
\begin{figure}
\centering
\begin{tikzpicture}
\begin{scope}[scale=2]
\draw[thick](0,0) circle (1cm);
\filldraw[fill=white, thick](.8623,0) circle (.5063cm);
\fill[white](0,0) circle (.99cm);
\draw[thick](0.8623,0.5063) -- (1.8623, 0.5063);
\draw[thick](0.8623,-0.5063) -- (1.8623, -0.5063);
\draw[dash pattern=on 1pt off 1pt] (0,0) -- (0.8623, 0);
\draw[dash pattern=on 1pt off 1pt] (0.8623, 0) -- (0.8623, 0.5063);
\draw[dash pattern=on 1pt off 1pt] (0,0) -- (0.8623, 0.5063);
\node at (0.45,0.4) {$1$};
\node at (1,0.2) {$r_0$};
\node at (0.4,0.1) {$\theta_0$};
\draw[dash pattern=on 1pt off 1pt] (0.25,0) arc (0:39:0.2);
\begin{scope}
\clip (1.1123,-0.6) rectangle (1.7,0.6);
\draw[dashed](1.1123,0) circle (.5063cm);
\end{scope}
\begin{scope}
\clip (1.3623,-0.6) rectangle (1.9,0.6);
\draw[dashed](1.3623,0) circle (.5063cm);
\end{scope}
\begin{scope}
\clip (1.6123,-0.6) rectangle (2.2,0.6);
\draw[dashed](1.6123,0) circle (.5063cm);
\end{scope}
\begin{scope}
\clip (1.8623,-0.6) rectangle (2.5,0.6);
\draw[thick](1.8623,0) circle (.5063cm);
\end{scope}
\end{scope}
\end{tikzpicture}
\caption{The Pinocchio set.}\label{fig:pinocchio}
\end{figure}

\section{A refined criterion}\label{sec:Saracco}

In this last section we state and prove an improved version of \cref{crit:Chen}. First, let us give a more general definition of strict interior rolling ball condition, that does not necessitate the piecewise Lipschitz regularity $\Omega$, rather it is stated through a lower bound on the reach of the complement set $\R^2 \setminus \Omega$, originally stated in~\cite[Def.~1.1]{Sar21}.
\begin{defin}
Let $\Omega\subset \R^2$ be a Jordan domain. We say that it enjoys the (weak) interior rolling ball condition of radius $r$ if $\RR(\R^2\setminus \Omega) \ge r$. We say that it enjoys the (weak) strict interior rolling ball condition if additionally for all $z \in \partial ((\R^2 \setminus \Omega) \oplus B_R)$ no antipodal points of $\partial B_R(z)$ lie both on $\partial \Omega$.
\end{defin}
For the sake of completeness, we recall that a Jordan curve is the image of a continuous and injective function $\Phi : \S^1 \to \R^2$ and a Jordan domain is the open region bounded by such a curve, and this is well defined thanks to the Jordan--Schoenflies Theorem. While any piecewise Lipschitz and simply connected set clearly is a Jordan domain, it might not be immediate to the reader unfamiliar with curvature measures that the condition on the reach is just a weaker request than that made previously on the existence of an interior rolling ball. The reach of a set $A$, first introduced in the foundational work~\cite{Fed59}, refer also to the survey~\cite{Tha08} and the comprehensive book~\cite{RZ19book}, is defined as follows. 

A set $A$ has reach $r$ if, for all $\rho <r$, the points in the \emph{Minkowski sum} $A\oplus B_\rho$ have a unique projection on $A$. Again, roughly speaking, this amounts to saying that it is possible to roll a ball of radius $r$ along $\partial A$ on the exterior of $A$, essentially providing a (weak) one-sided bound on its curvature. Notice the word \emph{exterior}: this is why the \emph{interior} rolling ball condition is defined through the reach of the complement set. We shall see that such a condition for $r=|\Omega|/P(\Omega)$, up to some very weak regularity condition on $\partial \Omega$, is sufficient for existence. Just as we did in the previous section, we notice that this weaker definition of interior rolling ball condition still implies that $\Omega$ has no necks of radius $r$, see~\cite[Lem.~3.1]{Sar21} where a finer result is proved. Assuming again this no neck condition, the criterion turns out to be necessary.
\begin{crit}\label{crit:via_reach}
Let $\Omega \subset \R^2$ be a Jordan domain such that 
\begin{equation}\label{eq:H(Om1)=0}
\calH^1(\Omega^{\scriptscriptstyle{(1)}} \cap \partial \Omega) = 0\,,
\end{equation} 
where $\Omega^{\scriptscriptstyle{(1)}}$ is the set of points of density $1$ for $\Omega$, and such that
\[
|\partial \Omega| = 0\,,
\]
where by this latter we mean that its boundary has zero \emph{$2$-dimensional Lebesgue} measure.

\begin{itemize}
\item[(i)] If it enjoys the (weak) strict interior rolling ball condition of radius $r=|\Omega|/P(\Omega)$, then the PDE~\eqref{eq:cap} with boundary condition~\eqref{eq:bc} has a solution for $\gamma =0$. 

\item[(ii)] If $\Omega$ has no necks of radius $r=|\Omega|/P(\Omega)$, then the PDE~\eqref{eq:cap} with boundary condition~\eqref{eq:bc} has a solution for $\gamma =0$ \emph{if and only if}  $\Omega$ enjoys the (weak) strict interior rolling ball condition of radius $r=|\Omega|/P(\Omega)$.
\end{itemize}
\end{crit}

\begin{proof}
By~\cite[Crit.~1.5]{Sar21} if $\Omega$ enjoys the (weak) strict interior rolling ball condition for $r=|\Omega|/P(\Omega)$, then $\Omega$ is the unique Cheeger set in itself, that is,
\[
h(\Omega) = \frac{P(\Omega)}{|\Omega|} < \frac{P(E)}{|E|}\,,
\]
for all proper subsets $E$. Hence, statement~(a) of \cref{thm:pmc-iff} holds. We are left with checking that the hypotheses of \cref{thm:pmc-iff} are met, as this would imply that statement~(c) of \cref{thm:pmc-iff} holds, which is our claim. We only need to check that~\eqref{eq:P=H}--\eqref{eq:Poincare} hold, as the other topological assumptions follow from $\Omega$ being a Jordan domain.

Since $\Omega$ is a Cheeger set and it satisfies~\eqref{eq:H(Om1)=0}, we can apply~\cite[Thm.~3.4]{Sar18} to find that $\Omega$ enjoys a Poincar\'e-type inequality, that is, hypothesis~\eqref{eq:Poincare} is met. Thus, applying~\cite[Lem.~3.5]{Sar18}, we find that
\[
\calH^1(\Omega^{\scriptscriptstyle{(0)}} \cap \partial \Omega) = 0\,.
\]
This latter equality, paired with~\eqref{eq:H(Om1)=0} and with the celebrated Federer's Structure Theorem implies that also hypothesis~\eqref{eq:P=H} is met.

To show the necessity when $\Omega$ has no necks of radius $r$, we reason as follows. If the PDE problem~\eqref{eq:cap}--\eqref{eq:bc} has a solution for $\gamma =0$, then the necessary condition
\[
\frac{P(\Omega)}{|\Omega|} < \frac{P(E)}{|E|}\,,
\]
holds, for all proper subsets $E$. Since Cheeger sets always exist when $\Omega$ is bounded, it remains proved that $\Omega$ is the unique Cheeger set in itself. Since we also have $|\partial \Omega|=0$, in virtue of~\cite[Thm.~1.4]{LNS17}, we have that $\Omega = \Omega^r \oplus B_r$, with $r=1/h(\Omega)$. The conclusion follows by~\cite[Lem.~3.1]{Sar21}.
\end{proof}

The regularity hypotheses requested in the three criteria went from piecewise Lipschitz to require that the set of points of density $1$ for $\Omega$ inside the boundary of $\Omega$ is negligible, paired with $|\partial \Omega|=0$. This latter request is of technical nature and due to the proof of~\cite[Thm.~1.4]{LNS17}. Even though it is not very stringent, as it forces us to discard only plane filling curves \`a la Knopp--Osgood (see~\cite[Chap.~8]{Sag94book}), we believe that it could be removed by employing the regularity theory of $\Lambda$-minimizers of the perimeter.

Finally, as mentioned at the beginning, we recall that existence of a solution of~\eqref{eq:cap}--\eqref{eq:bc} for the choice $\gamma=0$ also implies that for any choice $\gamma\in (0, \pi/2]$, which we sum up in the following corollary.

\begin{cor}
Let $\Omega \subset \R^2$ satisfy the general assumptions of~\cref{crit:via_reach} plus that of point~(i). Then, the PDE~\eqref{eq:cap}--\eqref{eq:bc} has solution for any choice $\gamma\in [0, \pi/2]$.
\end{cor}

In this case, we do not have (and we should not expect) an \emph{``if and only if''} statement. Indeed, it is well-known that there are sets $\Omega$ such that one has solutions of~\eqref{eq:cap}--\eqref{eq:bc} for all angles $\gamma\ge \gamma_0>0$, with $\gamma_0$ depending on the geometry of $\Omega$. This kind of phenomenon appears, e.g., when $\partial \Omega$ presents angles smaller than $\pi$. For the sake of completeness, we mention that there are some criteria taking into account the opening of the angle which provide existence for $\gamma\neq0$, refer to~\cite{Fin83, Tam86} and the examples therein.

%\section{Conclusions}
%

\section*{Acknowledgement}

G.S.\ is a member of INdAM and has been partially supported by the INdAM--GNAMPA Projects 2023 ``Esistenza e propriet\`a fini di forme ottime'' (codice CUP{\textunderscore}E53\-C22\-00\-19\-30\-001) and 2024 ``Ottimizzazione e disuguaglianze funzionali per problemi geometrico--spettrali locali e nonlocali'' (codice CUP{\textunderscore}E53\-C23\-00\-16\-70\-001).%
%
%\section*{Conflict of interests}
%
%\noindent The author declares no competing interest.

\bibliographystyle{plainurl}

\bibliography{capillarity_existence_geometric_criterion}

\end{document}